\documentclass[12pt,reqno,intlim]{amsart}

\usepackage{amssymb,amsmath}

\usepackage[arrow,matrix,curve]{xy}

\newdir{ >}{{}*!/-5pt/\dir{>}}

\newcommand{\nc}{\newcommand}

\nc{\bC}{\bold{C}}
\nc{\bN}{\Bbb{N}}
\nc{\cF}{\mathcal{F}}
\nc{\cE}{\mathcal{E}}
\nc{\cR}{\mathcal{R}}
\nc{\cM}{\mathcal{M}}
\nc{\al}{\alpha}
\nc{\bt}{\beta}
\nc{\gm}{\gamma}
\nc{\dl}{\delta}
\nc{\om}{\omega}
\nc{\sg}{\sigma}
\nc{\Sg}{\Sigma}
\nc{\vf}{\varphi}
\nc{\ve}{\varepsilon}
\nc{\os}{\overset}
\nc{\ol}{\overline}
\nc{\ul}{\underline}
\nc{\us}{\underset}
\nc{\sbs}{\subset}
\nc{\bsl}{\backslash}
\nc{\Ra}{\Rightarrow}
\nc{\lra}{\longrightarrow}
\nc{\all}{\allowdisplaybreaks}

\nc{\Codes}{\operatorname{{\bold{Codes}}}}
\nc{\RegMono}{\operatorname{\mathcal{R}{\rm{eg}\mathcal{M}{\rm{ono}\!}}}}
\nc{\RegEpi}{\operatorname{\mathcal{R}{\rm{eg}\mathcal{E}{\rm{pi}\!}}}}
\nc{\Mn}{\operatorname{\mathcal{M}{\rm{ono}\!}}}
\nc{\Ep}{\operatorname{\mathcal{E}{\rm{pi}\!}}}
\nc{\Rg}{\operatorname{\mathcal{R}{\rm{eg}\!}}}
\nc{\Ob}{\operatorname{Ob\!}}

\numberwithin{equation}{section}

\newtheorem{theo}{\ \ \ Theorem}[section]
\newtheorem{lem}[theo]{\ \ \ Lemma}
\newtheorem{prop}[theo]{\ \ \ Proposition}

\theoremstyle{definition}
\newtheorem{exmp}[theo]{\ \ \ Example}

\theoremstyle{remark}
\newtheorem{rem}[theo]{\ \ \ Remark}

\begin{document}

\title[]
{EFFECTIVE CODESCENT MORPHISMS  IN SOME VARIETIES OF UNIVERSAL
ALGEBRAS}

\author{Dali Zangurashvili}

\maketitle

\begin{abstract}
The paper gives a sufficient condition formulated in a syntactical
form for all codescent morphisms of a variety of universal algebras
satisfying the amalgamation property to be effective. This result is
further used in proving that all codescent morphisms of quasigroups
are effective.
\bigskip


\end{abstract}

\section{Introduction}

The present paper deals with the problem of describing effective
codescent morphisms in varieties of universal algebras. This problem
has been solved for modules over a commutative ring with unit (see
e. g. G. Janelidze and W. Tholen \cite{JT2}), commutative rings with
units (A. Joyal and M. Tierney (unpublished), B. Mesablishvili
\cite{M}), Boolean algebras (M. Makkai (unpublished)) and groups
\cite{Z1}. In this work we consider arbitrary varieties satisfying
the amalgamation property and give for them a sufficient condition
formulated in a syntactical form, for a codescent morphism to be
effective. In particular, this criterion implies the following:
\vskip+2mm

{\em Let a variety of universal algebras satisfy the amalgamation
property and the following condition:

$(*)$ Let
$$\xymatrix{B\ar[r]^-p\ar[d]_-q&E\ar[d]^{q'}\\
D\ar[r]^-{p'}&D\sqcup_BE}$$
 be a pushout with monomorphic $p$, $q$ (and $p'$,
$q'$). Then for any element $\alpha$ of $D\sqcup_{B}E$ and any
subalgebra $C$$^{1)}$\footnotetext{$^{1)}$In fact here, for any
$\alpha$, we can confine the consideration only to subalgebras $C$
which are generated by $B$ and elements $d_1,d_2,..., d_k$ from $D$,
for some presentation $\alpha=t(d_1, d_2,..., d_k, x_1, x_2, ...,
x_l)$ with $x_1, x_2, ..., x_l$ lying in $E$.} of $D$ containing $B$
and such that $\alpha$ lies in $C\sqcup_{B}E$
$^{2)}$\footnotetext{$^{2)}$One can easily show (see Lemma 2.3) that
$C\sqcup_{B}E$ is embedded to $D\sqcup_{B}E$.}, there exists a
presentation
$$\alpha=t(c_1, c_2,...,c_m, e_1, e_2,..., e_n)$$ in which
$^{3)}$\footnotetext{$^{3)}$We do not exclude the case where either
$m=0$ or $n=0$ (i.e., the corresponding variables are absent).} :
$c_1, c_2,..., c_m$ are in $C$; $e_1, e_2,..., e_n$ are in $E$; and
$t$ is an $(m+n)$-ary term such that, for any $d_1, d_2,..., d_m$ in
$D$, the equality
$$t(c_1, c_2,..., c_m,e_1, e_2,...,e_n)=t(d_1,
d_2,...,d_m,e_1,e_2,...e_n)$$ in $D\sqcup_{B}E$ implies that
$$d_1,d_2,...,d_m\in C.$$ \vskip+2mm

Then every codescent morphism of the variety is effective.}
 \vskip+5mm

The condition (*) is easily satisfied, for instance, for modules
over a commutative ring with unit and for groups. We show that it is
satisfied for quasigroups, too. Since the variety of quasigroups
satisfies also the amalgamation property, as proved by J. Je\v{z}ek
and T. Kepka in \cite{JK}, we obtain that

\;

{\em every codescent morphism of quasigroups is effective.}
\vskip+2mm

As to an internal characterization of such morphisms, we recall that
it is given in \cite{Z1} for any variety with the so-called strong
amalgamation property (and the variety of quasigroups is indeed of
this kind \cite{JK}).
 According to this
characterization, a monomorphism $p:B\rightarrowtail E$ is a
codescent morphism if and only if \;  $R'\cap(B\times B)=R$ \;
  for any congruence $R$ on $B$ and its closure $R'$
in $E$.

\section{Preliminaries}
We begin with the needed definitions from descent theory \cite{JT2}
formulated, for convenience, in the dual form.

Let $\bC$ be a category with pushouts, and let $p:B\to E$ be a
morphism in $\bC$. It is well known that the change-of-cobase
functor
$$p_*:B/\bC\lra E/\bC$$
(pushing out along $p$) has a right adjoint $p^!$ composing with $p$
from the right. $p$ is called a {\em codescent} (resp. {\em
effective codescent}) {\em morphism} if $p_*$ is precomonadic (resp.
comonadic), i.e., the comparison functor
$$\Phi_p:B/\bC\lra\Codes(p),$$
where $\Codes(p$) is the Eilenberg-Moore category of the comonad induced by
the adjunction
$$p_*\dashv \,p^!,$$
is full and faithful (resp. an equivalence of categories). Recall that
objects of $\Codes(p)$ (called codescent data with respect to $p$) are
triples $(C,\gm,\xi)$ with $C\in\Ob\bC$ and $\gm,\xi$ being morphisms
$E\lra C$ and \linebreak$C\lra C\sqcup_BE$, respectively, such that
the following equalities are valid (see Fig. 1 and 2):
\begin{align}
\xi\gm&=i_2,\\
(1_C,\gm)\xi&=1_C,\\
(i_1\sqcup_B1_E)\xi&=(\xi\sqcup_B1_E)\xi.
\end{align}

$$\xymatrix{B\ar[r]^-{p}\ar[d]_-{p}&E\ar[dd]^-{i_2}
\ar@/^1.73pc/[dddr]^-{\gm}\\
E\ar[d]_-{\gm}&\\
C\ar@/_1.4pc/[drr]_-{1_C}\ar@<0.4ex>[r]^-{i_1}
\ar@<-0.4ex>[r]_-{\xi}&C\sqcup_BE\ar[dr]^-{(1_C,\gm)}\\
&&C}\quad
\xymatrix{B\ar[r]^-{p}\ar[d]_-{p}&E\;\;\;\;\ar[dd]^-{i_2}
\ar@<1ex>[dddr]\\
E\ar[d]_-{\gm}&\\
C\ar@{-->}@<-0.4ex>[d]_-{i_1}\ar@<0.4ex>[d]^-\xi
\ar@<0.45ex>[r]^-{i_1}\ar@<-0.45ex>[r]_-{\xi}
&C\sqcup_BE\ar@{-->}@<-0.8ex>[dr]_{i_1\sqcup_B1_E}
\ar[dr]^(.4){\,\xi\,\sqcup_B1_E}\\
C\sqcup_BE\ar[rr]&&(C\sqcup_BE)\sqcup_BE}
$$
\centerline{\hskip-1.5cm Fig. 1\hskip+4.5cm Fig. 2}
\vskip+2mm

\noindent Moreover, a $\Codes(p)$-morphism
$(C,\gm,\xi)\to(C',\gm',\xi')$ is a $\bC$-morphism $h:C\lra C'$ such
that the following diagram commutes:
$$\xymatrix{&E\ar[dl]_-\gm\ar[dr]^-{\gm'}\\
C\ar[d]_-\xi\ar[rr]^-h& &C'\ar[d]^-{\xi'}\\
C\sqcup_BE\ar[rr]^-{h\sqcup_B1_E}&&C'\sqcup_BE}$$
\vskip+2mm

\noindent Also recall that the functor $\Phi_p$ maps every $f:B\lra D$ to
$$(D\sqcup_BE,{\ol{i}}_2,{\ol{i}}_1\sqcup_B1_E),$$
where ${\ol{i}}_1$ and ${\ol{i}}_2$ are the pushouts of $p$ and $f$,
respectively, along each other.

From now on it will be assumed that, in addition to pushouts, $\bC$
also has equalizers.

\begin{theo}[\cite{JT1}, \cite{JT2}]
$p$ is a codescent morphism if and only if it is a universal regular
monomorphism, i.e., a morphism such that any of its pushouts is a
regular monomorphism.
\end{theo}

Let $(C,\gm,\xi)$ be codescent data with respect to a morphism $p$. Consider
the equalizer
\begin{equation}
\xymatrix{Q\ar[r]^-q&C\ar@<0.45ex>[r]^-{i_1}\ar@<-0.45ex>[r]_-{\xi}
&C\sqcup_BE}
\end{equation}
It gives the equalizer in $B/\bC$:
\begin{equation}
\xymatrix{&B\ar@{-->}[dl]_-\dl\ar[d]^{\gm p}\ar[dr]^{i_2p}\\
Q\ar[r]_-q&C\ar@<0.45ex>[r]^-{i_1}\ar@<-0.45ex>[r]_-{\xi}
&C\sqcup_BE}
\end{equation}

\begin{theo}
For a codescent morphism $p$ and codescent data $(C,\gm,\xi)$ with respect
to $p$, the following conditions are equivalent:

{\rm(i)} $(C,\gm,\xi)\approx\Phi_p(f)$, for some $f\in\Ob B/\bC$;

{\rm(ii)} $(C,\gm,\xi)\approx\Phi_p(\dl)$;

{\rm(iii)} the functor $p_*$ preserves the equalizer $(2.5)$;

{\rm(iv)} the morphism $(q,\gm):Q\sqcup_B E\lra C$ is an isomorphism;

{\rm(v)} $q\sqcup_B1_E$ is a monomorphism and there exists a morphism
\\ $\theta:C\lra Q\sqcup_BE$ $($see Fig. $3)$ with
$$\xi=\big(q\sqcup_B1_E\big)\theta.$$
$$\xymatrix{B\;\ar@{ >-{>}}@<-0.2ex>[d]_-p\ar[dr]^-\dl\ar@{>-{>}}[rr]^-p&&E
\ar[dr]^{{\ol{i}}_2}\ar[dd]^(.3){i_2}&\\
E\ar[d]_-\gm&Q\ar[dl]_-q\ar[rr]^(.3){{\ol{i}}_1}&&Q\sqcup_BE\ar[dl]^-{q\sqcup_B1_E}\\
C\ar@{-->}@<0.4ex>[urrr]^(.4)\theta\ar@<0.45ex>[rr]^-(.6){\xi}
\ar@<-0.45ex>[rr]_{i_1}&&C\sqcup_BE&}$$
\centerline{Fig. $3$}
\vskip+2mm

If $p_*$ maps regular monomorphisms to morphisms $($of the coslice category
$B/\bC)$ with monomorphic underlying morphisms, then one can omit
``$q\sqcup_B1_E$ is a monomorphism and'' from {\rm (v)}.
\end{theo}

\begin{proof}The equivalence of the  conditions (i)--(iv) is well
known (see e.g. G. Janelidze and W. Tholen \cite{JT2}).
(iii)$\Ra$(v) follows from (2.3).

(v)$\Ra$(iv): We will show that the morphism $\theta$ is the inverse of $(q,\gm)$.
First we observe that
\begin{equation}
\theta(q,\gm){\ol{i}}_1=\theta\,q={\ol{i}}_1,
\end{equation}
as follows from the equalities
$$\big(q\sqcup_B1_E\big)\theta\,q=\xi\,q=i_1\,q=\big(q\sqcup_B1_E\big){\ol{i}}_1
.$$
Moreover,
\begin{equation}
\theta(q,\gm){\ol{i}}_2=\theta\,\gm={\ol{i}}_2,
\end{equation}
since by (2.1) we have
$$\big(q\sqcup_B1_E\big)\theta\,\gm=\xi\,\gm=i_2=
\big(q\sqcup_B1_E\big){\ol{i}}_2.$$ From (2.6) and (2.7) we obtain
that $\theta(q,\gm)=1_{Q\sqcup_BE}$. The equality
$(q,\gm)\theta=1_C$ follows from (2.2) and the trivial observation
that $(q,\gm)=(1_C,\gm)(q\sqcup_B1_E)$.
\end{proof}

Let $\cM$ be a class of $\bC$-morphisms. Recall (see e. g. E. W.
Kiss, L. M\'arki, P. Pr\"ohle and  W. Tholen \cite{KMPT}) that $\bC$
is said to satisfy the {\em amalgamation property} with respect to
$\cM$ if, for any pushout
$$\xymatrix{C\ar[r]^-\nu\ar[d]_-\mu&D\ar[d]^{\mu'}\\
C'\ar[r]^-{\nu'}&D'}$$ with $\mu,\nu\in\cM$, we have
$\mu',\nu'\in\cM$, too. If, in addition, any such pushout  is also a
pullback square (or, equivalently, in the case of varieties of
universal algebras, $\nu'(C')\cap \mu'(D)=\nu'\,\mu(C)$), then $\bC$
is said to satisfy the strong amalgamation property. If the class
$\cM$ is not indicated, it is meant to be the class of all
monomorphisms.

\begin{lem}
Let $\bC$ satisfy the amalgamation property with respect to a morphism class $\cM$ containing all regular monomorphisms. If $p$ is a codescent morphism, then the functor $p_*$ preserves $\cM$-morphisms $($in the obvious sense$)$.
\end{lem}

\begin{proof}
It is sufficient to observe that if the back square and the upper inclined one
in the diagram
$$\xymatrix{{B\ \ \ }\ar[r]^-p\ar@<-1ex>[dd]\ar[dr]&{\ \ \ \ \ E}\ar@<1ex>[dr]\ar@{-}[d]<2ex>
&\\
&{C'\ }\ar@<2ex>[d]\ar[r]\ar[dl]_-\mu
&C'\sqcup_BE\ar[dl]^{\mu\sqcup_B1_E}\\
{C\ \ \ }\ar[r]&{\ C\sqcup_BE}&}$$
\vskip+2mm
\noindent are pushouts, then so is the lower inclined square.
\end{proof}

Before continuing our discussion, let us recall the following recent result.

\begin{theo}[{\cite{Z2}}]
Let $(\cE,\cM)$ be a factorization system on $\bC$ with
$\cE\subset\Ep\bC$. A codescent morphism $p$ is effective if and
only if, for any morphism $p'$ lying in $\cM$ and being the pushout
of $p$ along an $\cE$-morphism and any codescent data
$(C',\gm',\xi')$ with respect to $p'$ such that $\gm'\,p'\in\cM$,
there exists $f'$ from the corresponding coslice category, such that
$(C',\gm',\xi')$ is isomorphic to $\Phi_{p'}(f')$. If \ $\bC$
satisfies the amalgamation property with respect to $\cM$, then the
statement remains valid provided that ``$\gm'\,p'\in\cM$'' is
replaced by ``$\gm'\in\cM$''.
\end{theo}

Similarly to (2.4), for codescent data $(C',\gm',\xi')$ with respect to
a morphism $p'$, let a monomorphism $q':Q'\lra C'$ be the equalizer of
the pair $(i'_1,\xi')$ with obvious $i'_1$. Applying Theorem 2.2, Lemma 2.3
and Theorem 2.4, we obtain

\begin{prop}
Let $(\cE,\cM)$ be a proper factorization system on $\bC$ $($i.e. a
factorization system with  $\cE\sbs\Ep\bC$ and $\cM\sbs\Mn\bC)$, and
let $\bC$ satisfy the amalgamation property with respect to $\cM$. A
codescent morphism $p$ is effective if and only if,
 for any morphism $p':B'\lra E'$ lying
in $\cM$ and being the pushout of $p$ along an $\cE$-morphism and
any codescent data $(C',\gm',\xi')$ with respect to $p'$ such that
$\gm'\in\cM$, there exists a morphism $\theta'$ with
$\xi'=\big(q'\sqcup_B1_{E'}\big)\theta'$ ${}^{4)}$
\footnotetext{${}^{4)}$ In the same manner one can derive criteria
for the effectiveness of $p$ by using also the conditions (ii) and
(iii) of Theorem 2.2, but for our further purposes we need only the
present criterion.}
\end{prop}

\;

For the {\em proof} we only observe that the inclusion $\cE\sbs\Ep\bC$
obviously implies the inclusion $\RegMono\bC\sbs\cM$.\hfill\qed

\section{Effective Codescent Morphisms in Varieties with the Amalgamation Property}

Throughout this section, if it is not specified otherwise, we assume that
$\bC$ is a variety of universal algebras (of any type $\cF$) satisfying
the amalgamation property.

Let $p:B\lra E$ be a codescent morphism in $\bC$, and let $R$ be a congruence
on $B$. It is well known that the pushout of $p$ along the projection
$B\lra B'=B/R$ is the obvious morphism $p':B'\lra E'=E/R'$, where $R'$ is
the closure of $R$ in $E$. It is obvious that $p'$ is a codescent morphism
as well.

Let $C'$ be an extension of $B'$. Below we will deal with the free product of
the algebras $C'$ and $E'$ with the amalgamated subalgebra $B'$ (see Fig.~4).
$$\xymatrix{B'\ar[r]^-{p'}\ar[d]&E'\ar[d]^{i'_2}\\
C'\ar[r]^-{i'_1}&C'\sqcup_{B'}E'}$$ \centerline{Fig. 4} \vskip+2mm

\noindent For convenience, when no confusion might arise, we
identify $C'$ and $E'$ with their images under $i'_1$ and $i'_2$,
respectively, and hence consider them as subalgebras of
$C'\sqcup_{B'}E'$ (whose intersection is, in general, wider than
$B'$). We adopt the similar convention for the free product of
several algebras with an amalgamated subalgebra.

All the terms considered below are meant to be of type $\cF$. When
using the notation $t(c_1,c_2,\dots,c_m,e_1,e_2,\dots,e_n)$ for a
term we do not exclude the case where either $m=0$ or $n=0$ (i.e.,
the corresponding variables are absent).

\begin{prop}
The following conditions are equivalent:

{\rm(i)} $p$ is effective;

{\rm(ii)} for any congruence $R$ on $B$ and any codescent data
$(C',\gm',\xi')$ (with respect to $p'$) with  monomorphic $\gm'$,
there exists a homomorphism \linebreak$\theta':C'\lra
Q'\sqcup_{B'}E'$ with
$$\xi'=\big(q'\sqcup_{B'}1_{E'}\big)\theta';$$

{\rm(iii)} for any congruence $R$ on $B$ and any codescent data
$(C',\gm',\xi')$ (with respect to $p'$) with monomorphic $\gm'$, one
has
$$\xi'(C')\sbs\big(q'\sqcup_{B'}1_{E'}\big)
\big(Q'\sqcup_{B'}E'\big);$$

{\rm(iv)} for any congruence $R$ on $B$, any codescent data
$(C',\gm',\xi')$ (with respect to $p'$) with monomorphic $\gm'$ and
any $c\in C'$, there exists a presentation of $\xi'(c)$ as a term
$t$ over the set $C'\cup E'$ of variables, such that, for any
variable $c'$ from $C'\bsl E'$ involved in $t$,   one has
$$\xi'(c')=c'.$$
\end{prop}

\begin{proof}
The equivalence (i)$\Leftrightarrow$(ii) is precisely the contents of Proposition 2.5 for $\cM=\Mn\bC$ (and $\cE=\RegEpi\bC$). (ii)$\Leftrightarrow$(iii) follows from Lemma 2.3, while (iii)$\Leftrightarrow$(iv) is obvious.
\end{proof}

\begin{prop}
For the following conditions, one has
{\rm(i)}$\Ra${\rm(ii)}$\Ra${\rm(iii)}.

{\rm(i)} For any congruence $R$ on $B$, any extensions
$B'\rightarrowtail C'\rightarrowtail D' $ and any element $\al$ of
$C'\sqcup_{B'}E'$ there exists a presentation
\begin{gather}
\alpha=t(c_1,c_2,\dots,c_m,e_1,e_2,\dots,e_n)
\end{gather}
 in which: $c_1$,
$c_2$, ..., $c_m$ are in $C'$; $e_1$, $e_2$, ..., $e_n$ are in $E'$;
and $t$ is an $(m+n)$-ary term such that, for any
$d_1,d_2,\dots,d_m$ in $D'$, the equality
\begin{gather}
t(c_1,c_2,\dots,c_m,e_1,e_2,\dots,e_n)=\notag\\
\quad\quad=t(d_1,d_2,\dots,d_m,e_1,e_2,\dots,e_n)
\end{gather}
in $D'\sqcup_{B'} E'$ implies that
\begin{gather}
d_1, d_2,..., d_m\in C';
\end{gather}

{\rm(ii)} for any congruence $R$ on $B$, any extension $C'$ of $B'$
and any element $\al$ of $C'\sqcup_{B'}E'$, there exists a
presentation $(3.1)$ in which: $c_1$, $c_2$, ..., $c_m$ are in $C'$;
$e_1$, $e_2$, ..., $e_n$ are in $E'$; and $t$ is an $(m+n)$-ary term
such that for any elements $d_1$, $d_2$, $\dots$, $d_m$ in
$C'\sqcup_{B'}E'$, the equality $(3.2)$ in
$C'\sqcup_{B'}E'\sqcup_{B'}E'$, where the variables
$e_1,e_2,\dots,e_n$ (in the case of their existence) in all their
occurrences on both sides of $(3.2)$ are considered as
representatives of the third cofactor of
$C'\sqcup_{B'}E'\sqcup_{B'}E'$ implies (3.3);

{\rm (iii)} $p$ is an effective codescent morphism.

\end{prop}

\begin{proof}
For (i)$\Ra$(ii) it is sufficient to take $D'=C'\sqcup_{B'}E'$.

(ii)$\Ra$(iii): We will verify the validity of the condition (iv) of
Proposition 3.1. To this end, consider any codescent data
$(C',\gm',\xi')$ with respect to $p'$ such that $\gm'$ is a
monomorphism, and any $c\in C'$. If $\xi'(c)\in C'$, then, according
to (2.2), $\xi'(c)=c$. If $\xi'(c)=e\in E'$, then the term $e$ is
obviously the desired one. Suppose that $\xi'(c)$ does not lie in
$C'\cup E'$. Consider the presentation (3.1) of $\al= \xi'(c)$. From
(2.3) we have the equality (3.2), where
$d_1=\xi'(c_1),d_2=\xi'(c_2),\dots,d_m=\xi'(c_m)$ and the variables
$e_1,e_2,\dots,e_n$ are considered as representatives of the second
$E'$ in $C'\sqcup_{B'}E'\sqcup_{B'}E'$. Therefore, for all $i$
$(1\leq i\leq m)$, we have $\xi'(c_i)\in C'$ and hence, by (3.2) we
have $\xi'(c_i)=c_i$.
\end{proof}

Proposition 3.2 immediately implies
\begin{theo}
Let $\bC$ be a variety of universal algebras that satisfies the
amalgamation property and the following condition:

$(*)$ Let
$$\xymatrix{B\ar[r]^-p\ar[d]_-q&E\ar[d]^{q'}\\
D\ar[r]^-{p'}&D\sqcup_BE}$$
 be a pushout in $\bC$ with monomorphic $p$, $q$ (and $p'$,
$q'$). Then for any element $\alpha$ of $D\sqcup_{B}E$ and any
subalgebra $C$ of $D$ containing $B$ and such that $\alpha$ lies in
$C\sqcup_{B}E$ (contained in $D\sqcup_{B}E$), there exists a
presentation $(3.1)$ of $\al$ in which: $c_1, c_2,..., c_m$ are in
$C$; $e_1, e_2,..., e_n$ are in $E$; and $t$ is an $(m+n)$-ary term
such that, for any $d_1, d_2,..., d_m$ in $D$, the equality $(3.2)$
in $D\sqcup_{B}E$ implies that
$$d_1,d_2,...d_m\in C.$$

Then all codescent morphisms of $\bC$ are effective.

\end{theo}

\begin{rem}
In (*), for any $\alpha$ from $D\sqcup_{B}E$, we can confine the
consideration only to subalgebras $C$ which are generated by $B$ and
elements $d_1,d_2,..., d_k$ from $D$, for some presentation
$\alpha=t(d_1, d_2,..., d_k, x_1, x_2, ..., x_l)$ with $x_1, x_2,
..., x_l$ lying in $E$.
\end{rem}

\begin{exmp}

(i) It is well-known that all codescent morphisms of the variety of
Abelian groups (more generally, modules over a commutative ring with
unit) are effective. The condition (*) is obviously satisfied in
this  case. \vskip+3mm

 (ii) In \cite{Z1} we have
shown that every codescent morphism of the variety of groups is
effective. Let us show that this variety satisfies the condition
(*). To this end, let us first recall the well-known fact related to
the free product $G$ of groups $G_1$ and $G_2$ with an amalgamated
subgroup $B$ (we assume that $G_1 \bigcap G_2=B$)\cite{K}.

For any right coset of $G_1$ and of $G_2$  by $B$, except for $B$
itself, we choose a representative. We denote the set of all chosen
representatives by $\textsl{A}$. Then every element of $G$ can be
uniquely written as a product
\begin{gather}
b\,  a_1\, a_2\,...\,a_n,
\end{gather}
where $n\geq0$, $b\in B$, all $a_j$ lie in $\textsl{A}$ and no two
$a_j, a_{j+1}$ belong to one and the same $G_j$. Form (3.4) is
called $\textsl{A}$-canonical. The procedure how an element
\begin{gather}
a'_1\, a'_2\,...\,a'_n,
\end{gather}
of $G$ (taken in uncancellable form) can be reduced to the canonical
form is described in \cite{K} (see also \cite{Z1}). Roughly
speaking, we, beginning from the right, pick out the left
$B$-coefficient from a current factor in (3.5), and then multiply it
to the left neighbor.

Let us now take $G_1=D$, $G_2=E$ and choose a set $\textsl{A}$ of
representatives of right cosets. Then take $G_1=C$, $G_2=E$ and
choose a set $\textsl{A}'$ of representatives such that
$\textsl{A}'$ $\subset$ $\textsl{A}$. In other words, for cosets of
both $C$ and $E$ by $B$ we take the already chosen representatives.
Let us consider any element $\alpha$ from $C\sqcup_{B} E$ and take
its $\textsl{A}'$-canonical form
\begin{gather}
b\, c_1\, e_1\, c_2\, e_2\,...\,e_{m-1}\,c_m,
\end{gather}
with $c_1, c_2,..., c_m\in C$ and $e_1, e_2,..., e_{m-1}\in E$. We
show that (3.6) is the desired representation of $\alpha$, where $b$
is considered as an element of $E$ (the proof of the statement in
the case where (3.6) is ended by an element of $E$ is similar). To
this end, let us consider $d_1, d_2,..., d_m\in D$ with
\begin{gather}
b\, c_1\, e_1\, c_2\, e_2\,...\,e_{m-1}\, c_m\, = \,b \,d_1\, e_1\,
d_2\, e_2\,...\,e_{m-1}\, d_m.
\end{gather}
Let us reduce the right hand part of (3.7) to the
$\textsl{A}$-canonical form. For that we consider the presentation
$$ d_m \,= \,b_m\, d'_m,$$
with $d'_m\in \textsl{A}$ and $b_m\in B$. From the uniqueness of the
canonical form, we conclude that $d'_m = c_m$ and hence $d_m \in C$.
Let $1\leq k<m$. We have
$$d_k\, b'_k = b_k\, d'_{k+1},$$
for some $d'_{k+1}\in \textsl{A}$ and $b_k, b'_k\in B$. Again, from
the uniqueness of the canonical form, we have $d'_{k+1} = c_{k+1}$,
which implies that $d_k\in C$.
\vskip+3mm

(iii) Let $\cF$ contain only nullary and unary operations. Then
\linebreak $D\sqcup_{B}E$, as a set, is obviously isomorphic to the
corresponding pushout $(D\sqcup_{B}E)_{\text{\bf Set}}$ in the
category of sets. Therefore, the condition (*) holds in that case,
too.
\vskip+3mm

(iv) Let $\bC$ be the variety of all algebras of type $\cF$. Then,
as is well known, each element of $D\sqcup_{B}E$ can be uniquely
presented as a term over $(D\sqcup_{B}E)_{\text{\bf Set}}$ such that
the variables of none of its subterm lie in one and the same
cofactor of $D\sqcup_{B}E$. This implies that the condition (*) of
Theorem 3.3 holds.
\vskip+3mm

 We conclude that all codescent
morphisms are effective in both cases (iii) and (iv).
\end{exmp}

\section{Effective Codescent Morphisms of Quasigroups}

Let us now pass to the case where $\bC$ is the variety of
quasigroups. Take its usual presentation $(\cF,\Sg)$. Recall that
here $\cF$ consists of three binary operations
$\text{\fontsize{8}{8pt}\selectfont$\circ$},/,\bsl$, while $\Sg$ is
the set of the identities
\begin{align}
\bsl\big(x_1,\text{\fontsize{8}{8pt}\selectfont$\circ$}(x_1,x_2)\big)=x_2,\\
/\big(\text{\fontsize{8}{8pt}\selectfont$\circ$}(x_1,x_2),x_2\big)=x_1,\\
\text{\fontsize{8}{8pt}\selectfont$\circ$}\big(x_1,\bsl(x_1,x_2)\big)=x_2,\\
\text{\fontsize{8}{8pt}\selectfont$\circ$}\big(/(x_1,x_2),x_2\big)=x_1.
\end{align}

Let $B$ and $A_i$ $(1\leq i\leq k,\;k\in \bN)$ be quasigroups such
that, for any $i$, $x_1$, $x_2\notin A_i$. For simplicity, it is
assumed that $A_i\cap A_j=B$, for any distinct $i,j$.

To distinguish between terms over $X=\{x_1,x_2\}$ and those over
$\os{k}{\us{i=1}{\cup}}A_i$, we will use different notations for
them: the former will be denoted by the capital letter $T$ (perhaps
with (co)indices), while the latters -- by the small-case letter $t$
(perhaps with (co)indices). Unless specified otherwise, we will use
the word ``term'' to mean a term over the set
$\os{k}{\us{i=1}{\cup}}A_i$.

In the set of terms we introduce transformations of the following
forms, which below will be called {\em reduction transformations}
or, for short, {\em reductions}:

{\rm(i)} if a term $t$ contains a subterm
\begin{equation}
f(a_1,a_2)
\end{equation}
with $f$ one of  $\text{\fontsize{8}{8pt}\selectfont$\circ$},/,\bsl$
and $a_1,a_2$ variables from one and the same $A_i$, then we replace
(4.5) in $t$ by the corresponding element of $A_i$; in that case the
term (4.5) is called the {\em replaced term} of the reduction, while
the value of (4.5) in $A_i$ is called the {\em replacing term} of
the reduction;

{\rm(ii)} if, for some identity
$$T=x_i$$
from (4.1)-(4.4), a term $t$ contains a subterm $t'$ obtained from
$T$ by replacing the variables $x_1$ and $x_2$ by some terms $t_1$
and $t_2$, then we replace $t'$ in $t$ by $t_i$. The subterm $t'$ is
called the {\em replaced term} of the reduction and $t_i$ is called
the {\em replacing term} of the reduction.

A reduction transformation is said to be performed on an occurrence  $o$ of an operation in $t$ if $o$ is the first (from the left) among
all occurrences of operations in the replaced term of the reduction.

A term is called {\em irreducible} if it admits no reduction
transformation.

\begin{lem}
For  any identity
$$T=x_i$$
from $(4.1)-(4.4)$, performing any sequence of reductions over the
term $t$ obtained from $T$ by replacing the variables $x_1$ and
$x_2$ by any irreducible terms $t_1$ and $t_2$, we arrive either at
a reducible term or at the term $t_i$.
\end{lem}

\begin{proof}
Let us assume that a term $t$ has, for instance, the form
\begin{equation}
\bsl\big(t_1,\text{\fontsize{8}{8pt}\selectfont$\circ$}(t_1,t_2)\big)
\end{equation}
with irreducible $t_1$ and $t_2$. If
$\text{\fontsize{8}{8pt}\selectfont$\circ$}(t_1,t_2)$, too, is
irreducible, then the only possible reduction is the one on the
depicted in (4.6) occurrence of the operation $\bsl$ and obviously
giving the term $t_2$. Let now
$\text{\fontsize{8}{8pt}\selectfont$\circ$}(t_1,t_2)$ be reducible.
Then it admits a reduction $r$ on the first occurrence of the
operation $\text{\fontsize{8}{8pt}\selectfont$\circ$}$. If $r$ is of
the form (i), then both $t_1$ and $t_2$ are variables from some
$A_i$ and the element of $A_i$ corresponding to $t$ is $t_2$. If $r$
is of the form (ii), then
$\text{\fontsize{8}{8pt}\selectfont$\circ$}(t_1,t_2)$ is equal
either to
$\text{\fontsize{8}{8pt}\selectfont$\circ$}(t_1\bsl(t_1,t_3))$ or to
$\text{\fontsize{8}{8pt}\selectfont$\circ$}(/(t_4,t_2),t_2)$. Let us
consider the first case. Then, performing the reduction $r$, from
$t$ we obtain the term $\bsl(t_1,t_3)$ which is irreducible and
equal to $t_2$. In a likewise manner we obtain $t_2$ in the second
case too.  The proof of the assertion in the case, where $t$ is
obtained from the left part of anyone of (4.2)--(4.4) by replacing
the variables $x_1$ and $x_2$ by some irreducible terms, is
analogous.
\end{proof}

Lemma 4.1 immediately implies

\begin{lem}
Any two reductions performed on one and the same occurrence $o$ of an operation in a term
give one and the same result  provided that $o$ has irreducible
arguments.
\end{lem}

\begin{lem}
Applying only reduction transformations, each term can be reduced to
a unique irreducible form.
\end{lem}

\begin{proof}
The existence of the needed term is obvious. In the particular case,
where a given term $t$ has the form described in Lemma 4.1, the
uniqueness immediately follows from this lemma. For the general case
we apply the principle of mathematical induction on the length $l$
of $t$.

If $l=1$, then the validity of the statement is clear. Assume that $l>1$ and
that the statement is valid for all terms whose length is less than $l$.
Let $t$ have the form
\begin{equation}
f(t_1,t_2),
\end{equation}
where $f$ is one of the operations
$\text{\fontsize{8}{8pt}\selectfont$\circ$},/,\bsl$, while $t_1$ and
$t_2$ are some terms. Let us introduce the following sequence $S$ of
reduction transformations: first both $t_1$ and $t_2$ are reduced
(in a certain manner) to irreducible forms and then the reduction is
performed, if possible, on the occurrence of $f$ depicted in (4.7).

Let $S'$ be any sequence of reductions applicable to (4.7) and
yielding an irreducible term. If the first transformation from $S'$
is the reduction of some  $t_i$, then, by the assumption of
induction, the transformations $S$ and $S'$ give one and the same
result.

Let the first member of $S'$ be the reduction on the first
occurrence of $f$ in $t$ . If this transformation is of the form
(i), then both $t_1$ and $t_2$ are merely variables from one and the
same $A_i$ and thus the results of $S$ and $S'$ coincide. Now assume
that the first member of $S'$ has the form (ii). Then $t$ is
obtained from the term $T$ by replacing the variables $x_1$ and
$x_2$ by some terms which without loss of generality can be assumed
to be irreducible. But Lemma 4.1 implies that both sequences $S$ and
$S'$ give one and the same term.
\end{proof}

\begin{lem}
Given two terms, if there are two reduction transformations such that
by applying one of them to one of these terms and the other transformation
to the second term we obtain one and the same term, then  there exists a term
that can be reduced by sequences of reduction transformations to both original terms.
\end{lem}

\begin{proof}
Let  reductions $r_1$ and $r_2$ transform terms $t_1$ and $t_2$ to a term $t$,
and let the replaced (resp. replacing) term of $r_i$ be $t_{i1}$ (resp.
$t_{i2}$), $i=1,2$. Assume that  $t_{12}$ and $t_{22}$ are not
subterms of each other in $t$. Then the desired term $t'$ is the term obtained
from $t$ by replacing $t_{12}$ by $t_{11}$ and $t_{22}$ by $t_{21}$.

Let $t_{12}$ be a subterm in $t$ of $t_{22}$. If $r_2$ is a
reduction transformation of the form (ii), then $t_{21}$ contains
both $t_{22}$ and $t_{12}$ as subterms. Hence, in that case the
desired term is obtained from $t_2$ by replacing $t_{12}$ in
$t_{21}$ by $t_{11}$.

Let $r_2$ be of the form (i). Then $t_{12}=t_{22}$. If $r_1$ is of
the form (ii), then the existence of the term we want to find
follows from the foregoing arguments. Otherwise it is obtained from
the term $t_1$ by replacing $t_{11}$ by
$$\hskip+4cm /\big(\text{\fontsize{8}{8pt}\selectfont$\circ$}(t_{11},\bsl(t_{12},
\text{\fontsize{8}{8pt}\selectfont$\circ$}(t_{21},t_{12}))),t_{12}\big).\hskip+3cm$$
\end{proof}

Let us now consider the free product of $A_i$ $(1\leq i\leq k)$ with
the amalgamated subquasigroup $B$. It is well known that it is
isomorphic to the quotient of the $\cF$-algebra of terms over the
set $\os{k}{\us{i=1}{\cup}}A_i$ of variables with respect to the
congruence $R$, where a term $t$ is $R$-equivalent to a term $t'$ if
and only if either $t=t'$ or $t'$ can be obtained from $t$ by a
sequence of transformations being reductions or their inverses.

\begin{lem}
For any two $R$-equivalent terms there exists a term which by
applying only reduction transformations can be reduced to both ones.
\end{lem}

The {\it proof}  easily follows from Lemma 4.4.\hfill\qed
\vskip+2mm

Lemma 4.3 and Lemma 4.5 immediately give rise to

\begin{lem}
Any element of the free product of quasigroups $A_i$ $(1\leq i\leq
k)$ with the amalgamated subquasigroup $B$ can be uniquely presented
as an irreducible term over the set $\os{k}{\us{i=1}{\cup}}A_i$ of
variables.
\end{lem}

Recall that the variety of quasigroups satisfies the amalgamation
property  (J. Je\v{z}ek and T. Kepka  \cite{JK}). Therefore, from
Lemma 4.6 and Theorem 3.3 we obtain

\begin{theo}
Every codescent morphism of quasigroups is effective.
\end{theo}

\begin{proof}
Let us verify the validity of the condition (*) of Theorem 3.3. For
any element $\al$ of $C\sqcup_{B}E$ consider its irreducible
presentation $t$ . Let the term $t'$ obtained from $t$ by replacing
all available  variables $c_1,c_2,\dots,c_m$ from $C$ by
$d_1,d_2,\dots,d_m$ from $D$ present $\al$ in $D\sqcup_{B}E$. If $t$
is irreducible, then, according to Lemma 4.6, $c_i=d_i$, for all
$i$. Assume now that $t'$ is not irreducible. Then reducing $t'$, we
obtain an irreducible term $t''$, $R$-equivalent to $t$ and such
that its length is less than that of $t$. This contradicts
Lemma~4.6.
\end{proof}

\section*{Acknowledgement}

 The author gratefully acknowledges  valuable discussions with George
Janelidze on the subject of this paper. The work was partially
supported by Volkswagen Foundation Ref.: I/84 328 and Georgian
National Science Foundation Ref.: ST06/3-004.

Authors address:
\vskip+2mm

Andrea Razmadze Mathematical Institute,

Tbilisi Centre for Mathematical Sciences

1 Alexidze Str., 0193, Tbilisi, Georgia

E-mail: dalizan@rmi.acnet.ge

\end{document}